\documentclass[a4paper]{article}

\usepackage{latexsym}
\usepackage{amsfonts}
\usepackage[reqno]{amsmath}
\usepackage{amsthm}
\usepackage{amssymb}

\usepackage{bbm}
\usepackage{mathrsfs}

\usepackage{pstricks}
\usepackage{pstricks-add}

\usepackage{ifthen}
\usepackage[nottoc]{tocbibind}

\usepackage[utf8]{inputenc}
\usepackage[T1]{fontenc}

\usepackage{textcomp}
\usepackage{lmodern}
\usepackage[a4paper]{geometry}

\usepackage[ps2pdf]{hyperref}
\hypersetup{colorlinks=true,linkcolor=blue}
  \usepackage[numbers,sort&compress]{natbib}

\newcounter{numtho}

\newtheoremstyle{mes_tho}%
		{9pt}{9pt}%
		{\itshape}%
		{}%
		{\bfseries}{.}%
		{\newline}%
		{}%

\theoremstyle{remark}

\theoremstyle{mes_tho}	\newtheorem{lemma}[numtho]{Lemma}
			\newtheorem{theorem}[numtho]{Theorem}
			
			\newtheorem{corollary}[numtho]{Corollary}
			\newtheorem{proposition}[numtho]{Proposition}

\newcommand{\arxiv}[2][]{%
\ifthenelse{\equal{#1}{}}{%
	\href{http://www.arxiv.org/abs/#2}{{\tt arXiv:#2}}%
			}
			{%
	\href{http://www.arxiv.org/abs/#2}{{\tt arXiv:#2}[#1]}%
}%
}				%

\newcommand{\defe}[2]{\textbf{#1}\index{#2}}
\newcommand{\tq}{\text{ st }}
\DeclareMathOperator{\SL}{SL}
\DeclareMathOperator{\SO}{SO}
\DeclareMathOperator{\so}{\mathfrak{so}}
\DeclareMathOperator{\ad}{ad}
\DeclareMathOperator{\Ad}{Ad}
\DeclareMathOperator{\AD}{\textbf{Ad}}
\DeclareMathOperator{\Adh}{Adh}
\DeclareMathOperator{\Int}{Int}

\newcommand{\hH}{\mathscr{H}}
\newcommand{\hS}{\mathscr{S}}			%
\newcommand{\hF}{\mathscr{F}}			%

\newcommand{\sA}{\mathcal{A}}
\newcommand{\sG}{\mathcal{G}}
\newcommand{\sH}{\mathcal{H}}			%
\newcommand{\sK}{\mathcal{K}}			%
\newcommand{\sN}{\mathcal{N}}
\newcommand{\sP}{\mathcal{P}}
\newcommand{\sQ}{\mathcal{Q}}
\newcommand{\sR}{\mathcal{R}}

\newcommand{\mO}{\mathcal{O}}

\newcommand{\eC}{\mathbbm{C}}

\newcommand{\eR}{\mathbbm{R}}

\newcommand{\mtu}{\mathbbm{1}}  			%

\makeindex

\begin{document}

\author{Laurent Claessens}
\title{The horizon of the BTZ black hole}
\maketitle

\begin{abstract}
	This paper is a sequel of \emph{Solvable symmetric black hole in anti de Sitter spaces} \cite{lcTNAdS}. In the latter, we described the BTZ black hole in every dimension by defining the singularity as the closed orbits of the Iwasawa subgroup of $\SO(2,n)$. In this article, we study the horizon of the black hole and we show that it is expressed as lateral classes of one point of the space. The computation is given in the four-dimensional case, but it makes no doubt that it can be generalized to any dimension.

	The main idea is to define an ``inclusion map'' from $AdS_3$ into $AdS_4$ and to show that all the relevant structure pass trough the inclusion. We prove, for example, that the inclusion of the three dimensional horizon into $AdS_4$ belongs to the four dimensional horizon : $\iota(\hH_3)\subseteq\hH_4$, then we deduce the expression of the horizon in $AdS_4$.
\end{abstract}

\tableofcontents

\section{Introduction}
\label{SecSumStructExist}

\subsection{Anti de Sitter space and the BTZ black hole}

The anti de Sitter space (hereafter abbreviated by $AdS$, or $AdS_l$ when we refer to a precise dimension) is a static solution to the Einstein's equations that describes a universe without mass. It was widely studied in different context in mathematics as well as in physics.

The BTZ black hole, initially introduced in \cite{BTZ_un,BTZ_deux} and then described and extended in various ways \cite{HolstPeldan,Aminneborg,Madden}, is an example of black hole structure which does not derives from a metric singularity. 

The structure of the BTZ black hole as we consider it here grown from the papers \cite{BTZB_deux,Keio} in the case of $AdS_3$. The dimensional generalization was first performed in \cite{lcTNAdS}. See also \cite{These} for for a longer review. Our point of view insists on the homogeneous space structure and the action of Iwasawa groups. One of the motivation in going that way is to embed the study of BTZ black hole into the noncommutative geometry and singleton physics \cite{BTZ_WZW,articleBVCS}. 

\subsection{The way we describe the BTZ black hole}

We look at the anti de Sitter space as the homogeneous space
\begin{equation}
	AdS_l=\frac{ G }{H }=\frac{ \SO(2,l-1) }{ \SO(1,l-1) }.
\end{equation}
We denote by $\sG=\so(2,l-1)$ and $\sH=\so(1,l-1)$ the Lie algebras and by $\pi$ the projection $G\to G/H$. The class of $g$ will be written $[g]$ or $\pi(g)$. We choose an involutive automorphism $\sigma\colon \sG\to \sG$ which fixes elements of $\sH$, and we call $\sQ$ the eigenspace of eigenvalue $-1$ of $\sigma$. Thus we have the reductive decomposition
\begin{equation}		\label{EqIntroRedDecompHQLieAlg}
	\sG=\sH\oplus\sQ.
\end{equation}
The compact part of $\SO(2,l-1)$ decomposes into $K=\SO(2)\times\SO(l-1)$.

Let $\theta$ be a Cartan involution which commutes with $\sigma$, and consider the corresponding Cartan decomposition
\begin{equation}
	\sG=\sK\oplus\sP,
\end{equation}
where $\sK$ is the $+1$ eigenspace of $\theta$ and $\sP$ is the $-1$ eigenspace. A maximal abelian algebra $\sA$ in $\sP$ has dimension two and one can choose a basis $\{ J_1,J_2 \}$ of $\sA$ in such a way that $J_1\in\sH$ and $J_2\in\sQ$.

Now we consider an Iwasawa decomposition
\begin{equation}
	\sG=\sK\oplus\sA\oplus\sN,
\end{equation}
and we denote by $\sR$ the Iwasawa component $\sR=\sA\oplus\sN$. We are also going to use the algebra $\bar\sN=\theta\sN$ and the corresponding Iwasawa component $\bar\sR=\sA\oplus\bar\sN$.

The Iwasawa groups $R=AN$ and $\bar R=A\bar N$ are naturally acting on anti de Sitter by $r[g]=[rg]$. It turns out that each of these two action has exactly two closed orbits, regardless to the dimension we are looking at. The first one is the orbit of the identity and the second one is the orbit of $[k_{\theta}]$ where $k_{\theta}$ is the element which generates the Cartan involution at the group level: $\AD(k_{\theta})=\theta$. In a suitable choice of matrix representation, the element $k_{\theta}$ is the block-diagonal element which has $-\mtu$ on $\SO(2)$ and $\mtu$ on $\SO(l-1)$. The $A\bar N$-orbits of $\mtu$ and $k_{\theta}$ are also closed. Moreover we have
\begin{equation}
	\begin{aligned}[]
		[A\bar N k_{\theta}]&=[k_{\theta}AN]\\
		[AN k_{\theta}]&=[k_{\theta}A\bar N]
	\end{aligned}
\end{equation}
because $A$ is invariant under $\AD(k_{\theta})$ and, by definition, $\Ad(k_{\theta})N=\bar N$. We define as \defe{singular}{Singular point} the points of the closed orbits of $AN$ and $A\bar N$ in $AdS$.

The Killing form of $\SO(2,l-1)$ induces a Lorentzian metric on $AdS$. The sign of the squared norm of a vector thus divides the vectors into three classes:
\begin{equation}
	\begin{aligned}[]
		\| X \|^2&>0&\rightarrow&&\text{time like,}\\
		\| X \|^2&<0&\rightarrow&&\text{space like,}\\
		\| X \|^2&=0&\rightarrow&&\text{light like.}
	\end{aligned}
\end{equation}
A geodesic is time (reps. space, light) like if its tangent vector is time like (reps. space, light).

If $E_1$ is a nilpotent element in $\sQ$, then every nilpotent in $\sQ$ are given by $\{ \Ad(k)E_1 \}_{k\in\SO(l-1)}$. These elements are also all the light like vectors at the base point. A light like geodesic trough the point $\pi(g)$ in the direction $\Ad(k)E_1$ is given by
\begin{equation}
	\pi(g e^{s\Ad(k)E_1}).
\end{equation}
One say that  points with $s>0$ are in the \defe{future}{Future} of $\pi(g)$ while points with $s<0$ are in the \defe{past}{Past} of $\pi(g)$. 

We say that a point in $AdS_l$ belongs to the \defe{black hole}{Black hole} if all the light like geodesics trough that point intersect the singularity in the future. We call \defe{horizon}{Horizon} the boundary of the set of points in the black hole. One say that there is a (non trivial) black hole structure when the horizon is non empty or, equivalently, when there are some points in the black hole, and some outside.

All these properties can be easily checked using the matrices given in \cite{These,lcTNAdS}. As far as notations are concerned, we denote by $X_{\alpha\beta}$ the basis of $\sN$ and $\bar\sN$ corresponding to our choice of Iwasawa decomposition. We have $\ad(J_1)X_{\alpha\beta}=\alpha X_{\alpha\beta}$ and $\ad(J_2)X_{\alpha\beta}=\beta X_{\alpha\beta}$.

\subsection{Organization of the paper}

In section \ref{SecOldResults}, we describe some old results about BTZ black hole.

In subsections \ref{SubSecProofExOld} and \ref{SubSecHorInThreeDimensionOld}, we recall how we proved the existence of the black hole structure in \cite{lcTNAdS} and how the horizon was described in the three dimensional case in \cite{Keio}. We adapt the latter result in our homogeneous space setting.

The subsection \ref{subSecTopoHor} gives some topological remarks about the black hole and the horizon. We point out that there are some light-like geodesics that are intersecting the singularity \emph{and then} the free part later in the future. We explain why that circumstance is very different from the situation of the most famous black holes in physics like the Schwarzschild's one.

Section \ref{SecNewWithMatrices} is devoted to the proof of our main result: the horizon of the BTZ black hole in $AdS_4$ is given by
\begin{equation}
	\hH_4=G_{X_{0+}}\cdot \iota(\hH_3)\cup G_{X_{0-}}\iota(\hH_3).
\end{equation}
where $\iota$ is the inclusion of $AdS_3$ in $AdS_4$ and $\hH_3$ is the horizon of the BTZ black hole in $AdS_3$.

\section{Some old results}
\label{SecOldResults}

\subsection{Proof of existence of a black hole}
\label{SubSecProofExOld}

Here is the way we proved in \cite{lcTNAdS} that the structure described in section \ref{SecSumStructExist} gives rise to a non trivial black hole. First, we see $AdS$ as embed in $\eR^{l+1}$ by the identification
\begin{equation}
	[g]=g e_u=g\begin{pmatrix}
		1	\\ 
		0	\\ 
		\vdots	\\ 
		0	
	\end{pmatrix}.
\end{equation}
If we name the coordinates as $(u,t,x,y,\ldots)$, we can prove that the singularity (closed orbits of $AN$ and $A\bar N$) is given by the equation $\hS\equiv t^2-y^2=0$. We can choose the matrices in such a way that nilpotent elements in $\sQ$ have the form
\begin{equation}		\label{EqAdkEMatr}
	E(w)=\begin{pmatrix}
		0	&	1	&	w_1	&	w_2	&	w_3	&	\ldots\\	
		-1	&		&		&		&		&	\\
		w_1	&		&		&		&		&	\\
		w_2	&		&		&		&		&	\\
		w_3	&		&		&		&		&	\\
		\vdots	&		&		&		&		&	\\			
	\end{pmatrix}
\end{equation}
with $\sum_iw_i^2=1$. Now, we can compute explicitly\footnote{by construction, $E(w)$ is nilpotent, so that the exponentiation is not a problem.} the value of $ e^{\mu q_0} e^{sE(w)}$ acting on $e_u$ and we write the expression $t(s)^2-y(s)^2$. Both $t(s)$ and $y(s)$ depend on the ``starting point'' $ e^{\mu q_0}$ and the direction $w$.

We prove that there exist some values of $\mu$ such that the solutions of $t(s)^2-y(s)^2=0$ are both positive for every $w$ (these points belong to the black hole), and we show that, for other values of $\mu$, we can find directions $w$ for which there are no solutions, or negative ones; the latter points are in the exterior of the black hole.

From that result, it is clear that a non trivial horizon exists. However, the question of the structure of the horizon was not yet addressed.

\subsection{Horizon in the three dimensional case}
\label{SubSecHorInThreeDimensionOld}

The structure of the horizon of $AdS_3$ was described in \cite{Keio} in the setting of $AdS_3=\SL(2,\eR)$. Our first job is to translate that result into the language of quotient of groups. This is done by the identification
\begin{equation}
	\begin{aligned}
		\psi\colon \SL(2,\eR)&\to AdS_3 \\
		\begin{pmatrix}
			u+x	&	y+t	\\ 
			y-t	&	u-x	
		\end{pmatrix}&\mapsto \begin{pmatrix}
			u	\\ 
			t	\\ 
			x	\\ 
			y	
		\end{pmatrix}.
	\end{aligned}
\end{equation}
We see that the points of the horizon are given by
\begin{equation}			\label{EqHOrAdSTroisVecteur}
	\begin{aligned}[]
		\pm\begin{pmatrix}
			\alpha	\\ 
			\cosh(a)	\\ 
			\alpha	\\ 
			\sinh(a)	
		\end{pmatrix}&&\text{and}&&\pm\begin{pmatrix}
			\alpha	\\ 
			\cosh(a)	\\ 
			-\alpha	\\ 
			\sinh(a)	
		\end{pmatrix},
	\end{aligned}
\end{equation}
which correspond to the points $(u,t,x,y)$ such that $u^2-x^2=0$. One should notice that these points can be expressed as lateral classes of the point $b=(0,1,0,0)$~:
\begin{equation}
	\hH_3=\pm G_{\{ J_1,X_{++} \}}b\cup\pm G_{\{ J_1,X_{--} \}}b
\end{equation}
where $G_{\{ X,Y \}}$ is the group of elements of the form $\exp(aX+bY)$. Notice that $G_{\{ J_1,X_{++} \}}b=G_{\{ J_1,X_{-+} \}}b$ and $G_{\{ J_1,X_{--} \}}b=G_{\{ J_1,X_{+-} \}}b$. For example,
\begin{equation}
	e^{aJ_2} e^{\alpha X_{++}}b=\begin{pmatrix}
		\alpha	\\ 
		\cosh(a)	\\ 
		\alpha	\\ 
		\sinh(a)	
	\end{pmatrix}.
\end{equation}
We are now intended to extend that result and express the horizon in $AdS_4$ as lateral classes of the horizon in $AdS_3$. Before to complete that work, we have to make a few remarks about the topology.

\subsection{Topology and horizon}
\label{subSecTopoHor}

The definition given in the previous sections produces a paradox. Let $x\in AdS$ and $l(s)$ be a light like geodesic trough $x$ which only intersects the singularity in past. We suppose that $l(0)=x$ and that $s_0<0$ is the biggest value of $s$ such that $l(s_0)\in \hS$. Thus, all points of the form $l(s)$ with $s_0<s<0$ are free. That form a sequence of free points which converges to the singularity, and then $l(s_0)$ belongs to the horizon.

This is however not possible in $AdS_3$ because the equation of the singularity is $t^2-y^2=0$ while the equation of the horizon is $u^2-x^2=0$. These two parts are really separated. 

The situation here is really different from the situation in the Schwarzschild's case. In the latter the singularity is well inside the horizon, and there are no geodesics reaching the infinity which have intersected the singularity in the past.

In our case, however, such geodesics do exist. The reason of such a difference resides in the fact that the causal structure (geodesics) are defined by the metric while, in our BTZ black hole, the singularity is not defined from metric considerations. There are thus no reasons to expect some compatibility relations like the fact to have a non naked singularity.

In order to correctly define the horizon, we have to introduce the space $BTZ=AdS\setminus\hS$ which in endowed with the induced topology. Then we define
\begin{equation}
	BH=\{ v\in BTZ\tq\forall k\in \SO(n),\, l_v^k(s)\in\hS\text{ has a solution with $s>0$} \}.
\end{equation}
Let us point out that the singularity itself is not part of the black hole, because it is not even part of $BTZ$. We define the free part of $BTZ$ as the set of points from which there exists a light-like geodesics which does not intersects the singularity in the future:
\begin{equation}
	\hF=\{ v\in BTZ\tq\exists k\in \SO(n),\, l_v^k(s)\in\hS\Rightarrow s<0 \}.
\end{equation}
The first definition makes that the black hole part is open by continuity and compactness of $\SO(n)$ : the minimum and the maximum of time to reach the singularity from one point of the black hole are both strictly positive numbers, and then can be maintained strictly positive in a neighborhood of the point.

\begin{proposition}		\label{PropBHouvertLibreFerme}
	The set of points in the black hole is open and set of free points is closed. In particular, the horizon is contained in the free set.
\end{proposition}

\begin{proof}
	The first point is the remark above. Now, the free part is closed in $BTZ$ as complementary of an open set.	
\end{proof}

The following theorem says that if the set of directions escaping the singularity from a point in $BTZ$ has an interior, then that point does not lies in the horizon. 
\begin{proposition}		\label{PropvFOsvghorvec}
	A point $v\in\hF_l$ such that there is an open set $\mO\subset S^{l-1}$ of directions for which $l^{w}_v(s)\in\hS$ has no solutions for $s\in\eR^+_0$ belongs to $\Int(\hF)$.
\end{proposition}

\begin{proof}
	Using the matricial representation \eqref{EqAdkEMatr}, we see that a point $v=[g]$ belongs to the singularity if the vector
	\begin{equation}
		g\cdot \begin{pmatrix}
			1	\\ 
			-s	\\ 
			s\bar w	
		\end{pmatrix}
	\end{equation}
	satisfies $t^2-y^2=0$. That equation is a second order polynomial in $s$ whose coefficients cannot be a constant for an open set with respect to $\bar w\in S^{l-1}$. From the assumptions, all the roots of that polynomial belong to $\eC\setminus\eR^+_0$. The latter being open, the roots of $l_{v'}^w(s)\in\hS$ are still in $\eC\setminus\eR^+_0$ when $v'$ runs over a small enough open set around $v$.

	We conclude that $v$ is in the interior of the free zone rather than on the horizon.
\end{proof}

An important characterisation of the horizon, pointed out in \cite{Keio}, is the following.
\begin{theorem}		\label{ThoHorIntDansS}
	A point belongs to the horizon if and only if the set of light-like directions for which the geodesics does not intersects the singularity has no interior in $S^{l-1}$.
\end{theorem}

\section{The horizon of the BTZ black hole}
\label{SecNewWithMatrices}

In this section, we show, that the horizon of the horizon of $AdS_4$ can be obtained using the action of a very simple group on the horizon of $AdS_3$, which is, itself, the orbit of one point under a known group. The result opens the possibility of describing the horizon in $AdS_l$ by induction on the dimension, and the possibility to compute the group which generates the horizon.  We define the inclusion map 
\begin{equation}
	\begin{aligned}
		\iota\colon AdS_3&\to AdS_4 \\
		\begin{pmatrix}
			u	\\ 
			t	\\ 
			x	\\ 
			y	
		\end{pmatrix}&\mapsto \begin{pmatrix}
			u	\\ 
			t	\\ 
			x	\\ 
			y	\\ 
			0	
		\end{pmatrix}.
	\end{aligned}
\end{equation}
At the matrix level, it corresponds to add a line and a column of zeros. We will denote by $\hF_l$ the free part of $AdS_l$. By definition, if $v\in\hF_l$, there exists a light like geodesic trough $v$ which does not intersect the singularity in the future. We also denote by $BH_l$ the set of elements of $AdS_l$ from which all the light-like geodesics intersect the singularity in the future.

Notice that $BH_l$ is open while $\hF_l$ is closed, as explained in proposition \ref{PropBHouvertLibreFerme}.

\begin{lemma}		\label{LemOouversttq}
	Let $v\in AdS_4$ and $g\in \SO(2,3)$ be a representative of $v$. If the set
	\begin{equation}
		\{ \begin{pmatrix}
			w_1	\\ 
			w_2	
		\end{pmatrix}\in S^2\tq
		\pi g\begin{pmatrix}
			1	\\ 
			-s	\\ 
			s\bar w	\\ 
			0	
		\end{pmatrix}\cap\hS_4=\emptyset\text{ with $s>0$}
				\}
	\end{equation}
	has an interior in $S^1$, then the set
	\begin{equation}
		\{ 
		\begin{pmatrix}
			w_1	\\ 
			w_2	\\ 
			w_3	
		\end{pmatrix}\in S^2\tq
		\pi g\begin{pmatrix}
			1	\\ 
			-s	\\ 
			s\bar w		
		\end{pmatrix}\cap\hS_4=\emptyset\text{ with $s>0$}
		\}
	\end{equation}
	has an interior in $S^2$.
\end{lemma}

\begin{proof}
The matrix $g$ in $\SO(2,3)$ representing the point $v$ has the form
\begin{equation}
	g=\begin{pmatrix}
 u	&	.	&	.	&	.	&	.\\ 
 t	&	a	&	b	&	c	&	d\\ 
 x	&	.	&	.	&	.	&	.\\ 
 y	&	a'	&	b'	&	c'	&	d'\\ 
z	&	.	&	.	&	.	&	. 
 \end{pmatrix}
\end{equation}
where the numbers $a,b,c,d,a',b',c',d'$ are not uniquely determined. We choose the representative in such a way to have $b\neq \pm b'$, which is always possible.

The assumption is that there exists an open set (with respect to $(w_1,w_2)\in S^1$) around $(w_1,w_2,0)$ such that the path
\begin{equation}		\label{EqPathgexpUTXYZ}
	\pi(g e^{s\Ad\left( k \right)E_1)})=
	\begin{pmatrix}
		U	\\ 
		T	\\ 
		X	\\ 
		Y	\\ 
		Z	
	\end{pmatrix}=
	\begin{pmatrix}
 u	&	.	&	.	&	.	&	.\\ 
 t	&	a	&	b	&	c	&	d\\ 
 x	&	.	&	.	&	.	&	.\\ 
 y	&	a'	&	b'	&	c'	&	d'\\ 
z	&	.	&	.	&	.	&	. 
 \end{pmatrix}
 \begin{pmatrix}
	 1	\\ 
	 -s	\\ 
	 sw_1	\\ 
	 sw_2	\\ 
	 0	
 \end{pmatrix}
\end{equation}
does not intersects the singularity in the future. In other words, we have $T\pm Y=0$ only with $s\leq 0$. Let
\begin{equation}
	\begin{aligned}[]
		T(w_1,w_2)&=t+s(bw_1+cw_2-a)\\
		Y(w_1,w_2)&=y+s(b'w_1+c'w_2-a')\\
		A_+(w_1,w_2)&=(b+b')w_1+(c+c')w_2-(a+a')\\
		A_-(w_1,w_2)&=(b-b')w_1+(c-c')w_2-(a-a').
	\end{aligned}
\end{equation}
We also denote by $\sigma_{\pm}$ the sign of $t\pm y$.

A simple computation shows that $T+Y=0$ when
\begin{equation}
	s=s_+=-\frac{ t+y }{ A_+(w_1,w_2) },
\end{equation}
and $T-Y=0$ when
\begin{equation}
	s=s_-=-\frac{ t-y }{ A_-(w_1,w_2) },
\end{equation}
The assumption is that the direction $(w_1,w_2,0)$ (and an open set in $S^1$ with respect to $(w_1,w_2)$) escapes the singularity, so that for every $(w_1',w_2')$ in a neighborhood of $(w_1,w_2)$, we have
\begin{equation}
	\begin{aligned}[]
		\sigma_{\pm}A_{\pm}(w_1',w_2')\geq 0,
	\end{aligned}
\end{equation}
which assures that the values of $s$ which annihilate $T+Y$ and $T-Y$ are negative or non existing. Since we choose $b\neq \pm b'$, the functions $A_{\pm}$ are nowhere constant, so we can find a direction $(w_1,w_2)$ such that $\sigma_{\pm}A_{\pm}(w_1,w_2)>0$. Notice that, by continuity, there exists a neighbourhood of $(w_1,w_2)$ in $S^1$ which escapes the singularity.

We are now studying what happens when one looks at a neighbourhood of $(w_1,w_2,0)$ in $S^3$. The path \eqref{EqPathgexpUTXYZ} is replaced by
\begin{equation}
	\pi(g e^{s\Ad(k)E_1})= 
	\begin{pmatrix}
 u	&	.	&	.	&	.	&	.\\ 
 t	&	a	&	b	&	c	&	d\\ 
 x	&	.	&	.	&	.	&	.\\ 
 y	&	a'	&	b'	&	c'	&	d'\\ 
z	&	.	&	.	&	.	&	. 
 \end{pmatrix}
\begin{pmatrix}
	1	\\ 
	-s	\\ 
	s(w_1+\epsilon_1)	\\ 
	s(w_2+\epsilon_2)	\\ 
	\epsilon_3	
\end{pmatrix},
\end{equation}
and we consider
\begin{equation}
	\begin{aligned}[]
		T(w_1,w_2,\bar\epsilon)&=t+s\big( b(w_1+\epsilon_1)+c(w_2+\epsilon_2)+d\epsilon_3-a \big)\\
		Y(w_1,w_2,\bar\epsilon)&=y+s\big( b'(w_1+\epsilon_1)+c'(w_2+\epsilon_2)+d'\epsilon_3-a' \big)
	\end{aligned}
\end{equation}
where $\bar\epsilon$ stands for $\epsilon_1$, $\epsilon_2$ and $\epsilon_3$. The same computations as before shows that $T+Y=0$ when
\begin{equation}
	s=s_+=-\frac{ t+y }{ A_+(w_1,w_2)+(b+b')\epsilon_1+(c+c')\epsilon_2+(d+d')\epsilon_3 },
\end{equation}
Since $\sigma_+A(w_1,w_2)>0$, there exists a $\delta$ such that $s_+$ remains negative for every choice of $\bar\epsilon<\delta$. The same holds with $T-Y$ which is zero when
\begin{equation}
	s=s_-=-\frac{ t-y }{ A_-(w_1,w_2)+(b-b')\epsilon_1+(c-c')\epsilon_2 +(d-d')\epsilon_3 }.
\end{equation}
Since $\sigma_-A_-(w_1,w_2)>0$, one can find a $\delta>0$ such that $\bar\epsilon<\delta$ implies that this fraction remains negative.

Thus, there exists a neighbourhood of $(w_1,w_2,0)$ in $S^2$ of directions escaping the singularity from the point $v$.
\end{proof}

\begin{lemma}		\label{LemIntTroisQueatr}
	With the notations defined before, we have
	\begin{equation}
		\iota\big( \Int(\hF_3) \big)\subseteq \Int\big( \hF_4 \big)
	\end{equation}
	where $\Int$ stands for the interior. In other words,
	\begin{equation}
		\Adh(BH_4)\cap\iota(AdS_3)\subset\iota\big( \Adh(BH_3) \big).
	\end{equation}
\end{lemma}

\begin{proof}

	Let $v=\iota(v')\notin\iota\big( \Adh(BH_3) \big)$, we also consider $g'$ a representative of $v'$ and $g=\iota(g')$, which is a representative of $v$. The element $v'$ is in the interior of the free zone: there exists an open set of directions which do not intersect the singularity of $AdS_3$ by theorem \ref{ThoHorIntDansS}. In other words, the set
\begin{equation}		\label{EqwwswswUn}
	\{ \begin{pmatrix}
	w_1	\\ 
	w_2	
\end{pmatrix}\in S^1\tq
\pi g'\begin{pmatrix}
	1	\\ 
	-s	\\ 
	sw_1	\\ 
	sw_2	
\end{pmatrix}\cap\hS_3 =\emptyset\}
\end{equation}
contains an open set of $S^1$. On the other hand, the $z$-component of the latter vector is obviously zero because $g=\iota(g')$ has the form
\begin{equation}
	g=\begin{pmatrix}
 .	&	.	&	.	&	.	&	0\\ 
 .	&	.	&	.	&	.	&	0\\ 
 .	&	.	&	.	&	.	&	0\\ 
 .	&	.	&	.	&	.	&	0\\ 
0	&	0	&	0	&	0	&	1 
 \end{pmatrix},
\end{equation}
thus equation \eqref{EqwwswswUn} can be ``extended'' and there exists an open set in $S^1$ such that
\begin{equation}
	\pi g\begin{pmatrix}
		1	\\ 
		-s	\\ 
		sw_1	\\ 
		sw_2	\\ 
		0	
	\end{pmatrix}\cap\iota(\hS_3)=\emptyset.
\end{equation}
Now, lemma \ref{LemOouversttq} shows that the set
\begin{equation}
	\{ 
		\begin{pmatrix}
			w_1	\\ 
			w_2	\\ 
			w_3	
		\end{pmatrix}\in S^2\tq
		\pi g\begin{pmatrix}
			1	\\ 
			-s	\\ 
			sw_1	\\ 
			sw_2	\\ 
			sw_3	
		\end{pmatrix}\cap\hS_4=\emptyset
	\}
\end{equation}
contains an open subset of $S^2$. That means that $\pi(g)=v$ belongs to the interior of $\hF_4$.
\end{proof}

\begin{proposition}		\label{PropFqTroisFt}
We have $\hF_4\cap\iota(AdS_3)\subset \iota(\hF_3)$.
\end{proposition}

\begin{proof}
Let $v\in\hF_4\cap\iota(AdS_3)$. With the same notations as above, we have
\begin{equation}		\label{EqRepresSOiotag}
	\iota(g')=
\begin{pmatrix}
 u	&	.	&	.	&	.	&	0\\ 
 t	&	a	&	b	&	c	&	0\\ 
 x	&	.	&	.	&	.	&	0\\ 
 y	&	a'	&	b'	&	c'	&	0\\ 
0	&	0	&	0	&	0	&	1 
 \end{pmatrix}
\end{equation}
The assumption is that, for every representative $g'$ of $v'$, there exists a direction $(w_1,w_2,w_3)\in S^2$ such that the path
\begin{equation}		\label{EqGedgpudt}
	\pi   \iota(g')\begin{pmatrix}
	1	\\ 
	-s	\\ 
	sw_1	\\ 
	sw_2	\\ 
	sw_3	
\end{pmatrix} 
\end{equation}
only intersects the singularity fore negative values of $s$. The values of $s$ that annihilate $t^2-y^2$ in the geodesic \eqref{EqGedgpudt} are
\begin{equation}
	\begin{aligned}[]
		s_+	&=-\frac{ t+y }{ -(a+a')+(b+b')w_1+(c+c')w_2 }\\
		s_-	&=-\frac{ t-y }{ -(a-a')+(b-b')w_1+(c-c')w_2 },
	\end{aligned}
\end{equation}
and these two values are either negative either non existing (vanishing denominator).

The work is now to find a direction $(w'_1,w'_2)\in S^1$ such that the geodesic
\begin{equation}
	\pi\big( g'\begin{pmatrix}
	1	\\ 
	-s	\\ 
	sw'_1	\\ 
	sw'_2	
\end{pmatrix} \big)
\end{equation}
does not intersect the singularity. The values of $s$ for which the latter geodesics intersects the singularity are
\begin{equation}
	\begin{aligned}[]
		s'_+	&=-\frac{ t+y }{ -(a+a')+(b+b')w'_1+(c+c')w'_2 }\\
		s'_-	&=-\frac{ t-y }{ -(a-a')+(b-b')w'_1+(c-c')w'_2 }.
	\end{aligned}
\end{equation}
If $w_3=0$, the proposition is true because one can choose $(w'_1,w'_2)=(w_1,w_2)$. If $w_3\neq 0$, the vector $(w_1,w_2)$ does not belong to $S^1$, and we have to find something else.

Let us consider the following two cases.
\begin{enumerate}
\item
there exists a representative \eqref{EqRepresSOiotag} with $a=a'=0$,
\item
there exists a representative \eqref{EqRepresSOiotag} with $c=c'=0$.
\end{enumerate}
In the first case, we have 
\begin{equation}		\label{EqDenoAAnnulerspm}
	s'_{\pm}=-\frac{ t\pm y }{ (b\pm b')w'_1+(c\pm c')w'_2 },
\end{equation}
and we can choose $(w'_1,w'_2)=N(w_1,w_2)$ with $N\in\eR$ fixed in such a way that $(w'_1,w'_2)\in S^1$. Thus we have $s'_{\pm}=\frac{1}{ N }s_{\pm}$ and it is sufficient to choose $N>0$ in order to leave the denominators of \eqref{EqDenoAAnnulerspm} of the right sign or zero.

In the second case, we have
\begin{equation}
	s'_{\pm}=-\frac{ t\pm y }{ -(a\pm a')+(b\pm b')w'_1 },
\end{equation}
thus one has to choose $w'_1=w_1$ and $w'_2=\sqrt{1-w_1^2}$.

Let us now discuss the values of $u$, $t$, $x$ and $y$ for which the first or the second cases are enforced. In order to be in the first case, we need to build a matrix of $\SO(2,2)$ of the form
\begin{equation}
	g'=\begin{pmatrix}
 u	&	\alpha	&	.	&	.	\\ 
 t	&	0	&	.	&	.	\\ 
 x	&	\beta	&	.	&	.	\\ 
 y	&	0	&	.	&	.	 
 \end{pmatrix}.
\end{equation}
That requires $\alpha^2-\beta^2=1$ and $u\alpha-x\beta=0$, while, for the second case, we need to build a matrix of $\SO(2,2)$ of the form
\begin{equation}
	g'=\begin{pmatrix}
 u	&	.	&	\alpha	&	.	\\ 
 t	&	.	&	0	&	.	\\ 
 x	&	.	&	\beta	&	.	\\ 
 y	&	.	&	0	&	.		 
 \end{pmatrix}.
\end{equation}
That requires $\alpha^2-\beta^2=-1$ and $u\alpha-x\beta=0$. 

In both cases, we have $\beta=\frac{ u }{ x }\alpha$ and $\alpha^2-\beta^2=\alpha^2\left( 1-\frac{ u^2 }{ x^2 } \right)$. If $| u |>| x |$, we can solve $\alpha^2-\beta^2=-1$, and if $| u |<| x |$, then we can solve $\alpha^2-\beta^2=1$. 

The last possible situation is $u=\pm x$. A point of $AdS_3$ in that situation belongs to the horizon by equation \eqref{EqHOrAdSTroisVecteur}, while one knows that point of horizon do have some directions which escape the singularity by corollary \ref{PropBHouvertLibreFerme}. Notice that in the latter situation, we do not use the assumption that $\iota(v')$ is free in $AdS_4$.
\end{proof}

\begin{corollary}		\label{CorBHBHHHHH}
	We have $\iota(BH_3)\subset BH_4$ and $\iota(\hH_3)\subset \hH_4$.
\end{corollary}

\begin{proof}
	If $\iota(v)\notin BH_4$, we have $\iota(v)\in \hF_4\cap\iota(AdS_3)\subset\iota(\hF_3)$, which is not possible if $v\in BH_3$.

	For the second part, we consider $v\in\hH_3\subset\hF_3$ (proposition \ref{PropBHouvertLibreFerme}). There is a direction $\begin{pmatrix}
		w_1	\\ 
		w_2	
	\end{pmatrix}\in S^1$ which escapes the singularity from $v$ in $AdS_3$. Of course, the direction $\begin{pmatrix}
		w_1	\\ 
		w_2	\\ 
		0
	\end{pmatrix}\in S^2$ escapes the singularity from $\iota(v)$ in $AdS_4$. Thus $\iota(v)\in\hF_4$.

	In every neighborhood of $v$, there exists a $\bar v\in BH_3$, and thus $\iota(\bar v)\in BH_4$. In other words, in every neighborhood of $\iota(v)$, there is that $\iota(\bar v)$ which belongs to $BH_4$. That proves that $\iota(v)$ belongs to $\hH_4$.
\end{proof}

\begin{lemma}
	We have $\hH_4\cap\iota(AdS_3)\subset\iota(\hH_3)$.
\end{lemma}

\begin{proof}
	Let $v\in\hH_4\cap\iota(AdS_3)$. Since $\hH_4\subset\hF_4$, we have $v\in\hF_4\cap\iota(AdS_3)\subset\iota(\hF_3)$ (proposition \ref{PropFqTroisFt}), and then there exists a $v'\in\hF_3$ such that $v=\iota(v')$. Now, we have to prove that $v'\in\hH_3$. If $v'$ belongs to the interior of $\hF_3$, lemma \ref{LemIntTroisQueatr} implies that
	\begin{equation}
		v=\iota(v')\in\iota\big( \Int(\hF_3) \big)\subset\Int(\hF_4),
	\end{equation}
	which disagrees with the fact that $v\in\hH_4$.
\end{proof}

\begin{proposition}		\label{PropovHhnonXYzero}
	Let $v'=(u',t',x',y',z')\in\hH_4$ with $u'$ and $x'$ not both vanishing. Then
	\begin{equation}
		v'\in G_{X_{0+}}\cdot \iota(\hH_3)\cup G_{X_{0-}}\cdot \iota(\hH_3).
	\end{equation}
\end{proposition}

\begin{proof}
	As a first step, we want to solve the equation
	\begin{equation}
		e^{\alpha X_{0+}}\begin{pmatrix}
			u	\\ 
			t	\\ 
			x	\\ 
			y	\\ 
			0	
		\end{pmatrix}=
		\begin{pmatrix}
			\frac{ \alpha^2(u-x) }{2}+u	\\ 
			t	\\ 
			\frac{ \alpha^2(u-x) }{2}+x	\\ 
			y	\\ 
			-\alpha(x-u)	
		\end{pmatrix}=\begin{pmatrix}
			u'	\\ 
			t'	\\ 
			x'	\\ 
			y'	\\ 
			z'	
		\end{pmatrix}
	\end{equation}
	with respect to $u$, $t$, $x$, $y$ and $\alpha$. The result is $t=t'$, $y=y'$ and
	\begin{equation}
		\begin{aligned}[]
			\alpha&=\frac{ z' }{ u'-x' },&u&=u'-\frac{ z'^2 }{ 2(u'-x') },&x&=\frac{ z'^2 }{ 2(u'-x') }-x'.
		\end{aligned}
	\end{equation}
	We conclude that, as long as $u'-x'\neq 0$, the point $v'$ belongs to $G_{X_{0+}}\cdot\iota(AdS_3)$. The same computation shows that $v'\in G_{X_{0-}}\cdot\iota(AdS_3)$ as long as $x'+u'\neq 0$. Let us observe that the actions of the matrices $ e^{\alpha X_{0+}}$ and $ e^{\beta X_{0-}}$ do not change the $t$ and $y$ component of a vector in $\eR^{2,l-1}$, so that the set of directions for which $v$ falls in the singularity is exactly the same as the set of directions for which $ e^{\alpha X_{0+}}v$ and $ e^{\beta X_{0-}}v$ fall in the singularity.
	
	Now, let us suppose that $v= e^{\alpha X_{0+}}\iota(v')$ with $v'\in AdS_3$. We want to prove that $\iota(v')\in\hH_4$ (i.e. there is an element in the black hole in each neighbourhood of $\iota(v')$) because in that case, corollary \ref{CorBHBHHHHH} would conclude that $v'\in\hH_3$.

	Let $\mO$ be a neighbourhood of $\iota(v')$. The set $ e^{\alpha X_{0+}}\mO$ is a neighborhood of $v$, and thus there exists an element $\bar v\in e^{\alpha X_{0+}}\mO\cap BH_4$. Now the element $ e^{-\alpha X_{0+}}\bar v$ belongs to $\mO\cap BH_4$, so that $\iota(v')$ belongs to $\hH_4$.
\end{proof}

\begin{lemma}		\label{Lemuxznonsing}
	A point of the form
	\begin{equation}
		v=\begin{pmatrix}
			0	\\ 
			t	\\ 
			0	\\ 
			y	\\ 
			z	
		\end{pmatrix}
	\end{equation}
	does not belongs to the horizon.
\end{lemma}

\begin{proof}

Since the horizon is $A$-invariant, we can reduce the lemma to the case of any element of the form $ e^{\eta J_1}v$. We have
\begin{equation}
	 e^{\eta J_1}
\begin{pmatrix}
	0	\\ 
	t	\\ 
	0	\\ 
	y	\\ 
	z	
\end{pmatrix}=
\begin{pmatrix}
 1	&	0		&	0	&	0		&	0\\ 
 0	&	\cosh(\eta)	&	0	&	\sinh(\eta)	&	0\\ 
 0	&	0		&	1	&	0		&	0\\ 
 0	&	\sinh(\eta)	&	0	&	\cosh(\eta)	&	0\\ 
 0	&	0		&	0	&	0		&	1 
 \end{pmatrix}
\begin{pmatrix}
	0	\\ 
	t	\\ 
	0	\\ 
	y	\\ 
	z	
\end{pmatrix}=
\begin{pmatrix}
	0				\\ 
	\cosh(\eta)t+\sinh(\eta)y	\\ 
	0				\\ 
	\sinh(\eta)t+\cosh(\eta)y	\\ 
	z
\end{pmatrix}
\end{equation}
We annihilate the $y$ component by choosing $\eta=\ln\left( \frac{ t-y }{ t+y } \right)$. The logarithm makes sense because, since $t^2-y^2-z^2=1$, we have $t^2-y^2\geq 0$. The case $t^2-y^2=0$ is trivial (the point $v$ belongs to the horizon), so that we can assume $t^2-y^2>0$.

A representative of $(0,t,0,0,z)$ in $\SO(2,2)$ is easy to find, and the geodesic in the direction $\bar w\in S^2$ is given by
\begin{equation}
	\begin{pmatrix}
 0	&	1	&	0	&	0	&	0\\ 
 t	&	0	&	0	&	0	&	-z\\ 
 0	&	0	&	1	&	0	&	0\\ 
 0	&	0	&	0	&	1	&	0\\ 
z	&	0	&	0	&	0	&	-t 
 \end{pmatrix}
\begin{pmatrix}
	1	\\ 
	-s	\\ 
	sw_1	\\ 
	sw_2	\\ 
	sw_3	
\end{pmatrix}=
\begin{pmatrix}
	.	\\ 
	t-szw_3	\\ 
	.	\\ 
	sw_2	\\ 
	.	
\end{pmatrix}.
\end{equation}
It belongs to the singularity when $s$ takes one of the values
\begin{equation}
	s_{\pm}=\frac{ t }{ w_3z\pm w_2 }.
\end{equation}
As long as $|w_2|<|w_3z|$, the two values $s_{\pm}$ have the same sign, which can be decided by making $w_3$ positive or negative. That provides an open set in $S^2$ of directions which escape the singularity, so that $v\notin\hH_4$ by proposition \ref{PropvFOsvghorvec}.
\end{proof}

\begin{theorem}			\label{ThoHorQuatreInclusionHorTrois}
	The horizon of $AdS_4$ is given by
	\begin{equation}		\label{EqEqHOrGVGXQuatr}
		\hH_4=G_{X_{0+}}\cdot \iota(\hH_3)\cup G_{X_{0-}}\iota(\hH_3).
	\end{equation}
	i.e. an union of lateral classes of the horizon of $AdS_3$ by one dimensional subgroups of $N$ and $\bar N$.

	The equation in the ambient $\eR^5$ is $\hH_4\equiv u^2-x^2-z^2=0$.
\end{theorem}

\begin{proof}
	We begin by the direct inclusion. If $v=(u,t,x,y,z)\in\hH_4$ with $u\neq 0$ or $x\neq 0$, we proved in proposition \ref{PropovHhnonXYzero} that $v$ has the form \eqref{EqEqHOrGVGXQuatr}. Now, if $u=x=0$, the lemma \ref{Lemuxznonsing} shows that $v$ does not belongs to the horizon.

	For the reverse inclusion, we know that elements of $\iota(\hH_3)$ belong to $\hH_4$ by corollary \ref{CorBHBHHHHH}. If $v$ belong to $\hH_4$, then $ e^{\alpha X_{0+}}v$ and $ e^{\beta X_{0-}}v$ also belong to the horizon.
\end{proof}

\section{Conclusion}

The horizon of the BTZ black hole in $AdS_3$ was already expressed in \cite{Keio} as lateral classes of one point under the action of the Iwasawa component of the isometry group of $AdS_3$.

We proved that the simple inclusion map $\iota\colon AdS_3\to AdS_4$ transports the causal structure (free zone, black hole, horizon) from $AdS_3$ to $AdS_4$. We studied in particular the way the horizon changes when ones jumps from dimension $3$ to dimension $4$ and we obtained that the horizon in $AdS_4$ is expressed as lateral classes of the inclusion of the horizon of $AdS_3$ in $AdS_4$. In the same time, we obtained a simple equation for the horizon seen as a subset of $\eR^5$.

Although the results are quite satisfying, the method used here to prove them is quite unsatisfactory because we didn't used all the wealth structure of $\so(2,3)$ and of its reductive decomposition \eqref{EqIntroRedDecompHQLieAlg}. We plan, in a future work, to get a much deeper understanding of the structure of $\sG$ and $\sQ$, in such a way to provide simpler proofs, in the same time as a dimensional generalization of the result of theorem \ref{ThoHorQuatreInclusionHorTrois}. We would also like to define a class of homogeneous spaces $G/H$ which accept a BTZ-like black hole.

\pagestyle{headings}

\bibliographystyle{unsrt}			

\begin{thebibliography}{10}

\bibitem{lcTNAdS}
Laurent Claessens and Stephane Detournay.
\newblock Solvable symmetric black hole in anti-de {S}itter spaces.
\newblock {\em J. Geom. Phys.}, 57:991--998, 2007.
\newblock \arxiv{math.DG/0510442}.

\bibitem{BTZ_un}
Máximo Ba{\~n}ados, Claudio Teitelboim, and Jorge Zanelli.
\newblock The lack hole in three-dimensional space-time.
\newblock {\em Phys. Rev. Lett.}, 69(13):1849--1851, 1992.
\newblock \arxiv{hep-th/9204099v3}.

\bibitem{BTZ_deux}
Máximo Ba{\~n}ados, Marc Henneaux, Claudio Teitelboim, and Jorge Zanelli.
\newblock Geometry of the {$2+1$} black hole.
\newblock {\em Phys. Rev. D (3)}, 48(4):1506--1525, 1993.
\newblock \arxiv{gr-qc/9302012v1}.

\bibitem{HolstPeldan}
S.~Holst and P.~Peld{\'a}n.
\newblock Black holes and causal structure in anti-de {S}itter isometric
  spacetimes.
\newblock {\em Classical Quantum Gravity}, 14(12):3433--3452, 1997.
\newblock \arxiv{gr-qc/9705067}.

\bibitem{Aminneborg}
S.~{\AA}minneborg, I.~Bengtsson, S.~Holst, and P.~Peld{\'a}n.
\newblock Making anti-de {S}itter black holes.
\newblock {\em Classical Quantum Gravity}, 13(10):2707--2714, 1996.
\newblock \arxiv{gr-qc/9604005}.

\bibitem{Madden}
O.~Madden and S.~F. Ross.
\newblock Quotients of anti-de {S}itter space.
\newblock {\em Phys. Rev. D (3)}, 70(2):026002, 8, 2004.
\newblock \arxiv{hep-th/0401205}.

\bibitem{BTZB_deux}
P.~Bieliavsky, S.~Detournay, M.~Herquet, M.~Rooman, and Ph. Spindel.
\newblock Global geometry of the $2+1$ rotating black-hole.
\newblock 2003.
\newblock \arxiv{hep-th/0306293v1}.

\bibitem{Keio}
Pierre Bieliavsky, Stephane Detournay, Philippe Spindel, and Marianne Rooman.
\newblock Noncommutative locally anti-de {S}itter black holes.
\newblock In {\em Noncommutative geometry and physics}, pages 17--33. World
  Sci. Publ., Hackensack, NJ, 2005.
\newblock \arxiv{math.QA/0507157}.

\bibitem{These}
Laurent Claessens.
\newblock Locally anti de {S}itter spaces and deformation quantization.
  {P}h.{D}. thesis.
\newblock September 2007.
\newblock \arxiv[math-DG]{0912.2215}.

\bibitem{BTZ_WZW}
Pierre Bieliavsky, Stephane Detournay, Marianne Rooman, and Philippe Spindel.
\newblock {BTZ} black holes, {WZW} models and noncommutative geometry.
\newblock 2005.
\newblock \arxiv{hep-th/0511080}.

\bibitem{articleBVCS}
Pierre Bieliavsky, Yannick Voglaire, Laurent Claessens, and Daniel Sternheimer.
\newblock Quantized anti de {S}itter spaces and non-formal deformation
  quantizations of symplectic symmetric spaces.
\newblock {\em Contemporary Mathematics}, (450), 2008.
\newblock \arxiv[math.QA]{0705.4179v1}.

\end{thebibliography}
\def\polhk#1{\setbox0=\hbox{#1}{\ooalign{\hidewidth
  \lower1.5ex\hbox{`}\hidewidth\crcr\unhbox0}}}

\end{document}